\documentclass[oneside,draft,letterpaper,12pt]{amsart}
\usepackage{amsthm}
\usepackage{amsmath}
\usepackage{amssymb}
\theoremstyle{definition}

\def\R{\text{$\mathbb{R}$}}

\def\fis#1{\dot{#1}}
\def\lra{\longrightarrow}

\def\fis#1{\dot{#1}}
\def\lra{\longrightarrow}

\def\d1#1#2{\frac{d#1}{d#2}}

\def\p1#1#2{\frac{\partial #1}{\partial #2}}

\def\to{t_o}
\def\pr{\parallel}

\def\part{a=\to \leq t_1 \leq \ldots \leq t_{n-1} \leq t_{n}=b}
\def\sg{\sum_{g \in G}}
\def\sgj{\sum_{g \in J}}

\newcommand{\dem}{{\bf Proof:}}
\newcommand{\n}{\noindent}

\def\C{\text{$\mathbb{C}$}}
\def\N{\text{$\mathbb{N}$}}
\def\E{\text{$\mathbb{E}$}}
\def\H{\text{$\mathbb{H}$}}
\def\Z{\text{$\mathbb{Z}$}}
\def\Q{\text{$\mathbb{Q}$}}

\def\vai{\rightarrow}

\title[Isometric actions]{Properly discontinuous 
isometric actions on the unit sphere of infinite dimensional Hilbert spaces}
\author[L. Biliotti]{Leonardo Biliotti}
\thanks{ Research partially supported by CNPq (Brazil)}
\begin{document}
\theoremstyle{plain}\newtheorem{thm}{Theorem}
\theoremstyle{plain}\newtheorem{prop}[thm]{Proposition}
\theoremstyle{plain}\newtheorem{lemma}[thm]{Lemma}
\theoremstyle{plain}\newtheorem{cor}[thm]{Corollary}
\theoremstyle{definition}\newtheorem{defini}[thm]{Definition}
\theoremstyle{remark}\newtheorem{remark}[thm]{Remark}
\theoremstyle{plain} \newtheorem{assum}[thm]{Assumption}
\theoremstyle{definition}\newtheorem{example}[thm]{Example}
\begin{abstract}
We study the properly discontinuous and isometric              
actions on the unit sphere of  infinite dimensional Hilbert spaces and   
we get some new examples of Hilbert manifold with constant 
positive sectional curvature. We prove 
some necessary conditions for a group to act isometrically 
and properly discontinuously, 
and in the case of finitely generated Abelian groups the necessary and 
sufficient conditions are given. 
\end{abstract}
\maketitle

\section{Introduction}
\mbox{}
A Hilbert manifold $M$ is a manifold modeled on a Hilbert space $\H$
and equipped with an inner product $\langle \ , \  \rangle (p),$
on each tangent space $T_p M \cong \H,$ 
depending smoothly on $p$ and defining
on $T_p M$ a norm equivalent to the original norm of $\H.$ 
In what follows we assume that $M$ is complete metric space with respect
to the usual distance obtained from Riemannian metric 
$ \langle \, \ \rangle.$  
In infinite dimensional geometry  most of the local results follow from 
general arguments analogous to those in the finite dimensional case 
(see \cite{Kl} or \cite{La}). The investigation of global properties,
on the contrary, is  
harder than the finite dimensional case; for example the  
Hopf-Rinow Theorem is only generically satisfied on a complete 
Hilbert manifold (see \cite{Ek}). An important 
application, due to Anderson (\cite{An}), is the extension of
the Bonnet Theorem 
about the estimate of the diameter of  complete Hilbert manifolds 
with sectional curvature $K \geq K_o >0$. In  finite dimensional geometry 
the Bonnet Theorem implies also that the fundamental group 
must be finite, (see \cite{Mdc}), while in \cite{An} doesn't appear any 
information about the fundamental group. 
As in the finite dimensional case there is a 
bijective correspondence between complete Hilbert manifolds, modeled on
$\H$, with positive constant 
sectional curvature $1$ 
and groups $G$ acting freely, isometrically and properly 
discontinuously on the unit sphere $\bf{S} (\H \times \R)$ in 
$\H \times \R.$  
We recall that a group $G$ 
acts {\em properly discontinuously} on a topological space $X$ if it 
satisfies the following conditions:
\begin{enumerate}
\item if $x,y \in X,$ $ y \not\in Gx,$ then $x$ and $y$ have 
neighborhoods $U$ and $V$ such that $gU \cap V= \emptyset$ for every 
$g \in G$; 
\item for any $x \in X$ the set  $G_x:=\{ g \in G:\ gx=x \}$  is finite;  
\item for any $x \in X$  there exists a neighborhood $U$, stable by $G_x$,
such that $gU \cap U= \emptyset$ for every $g \in G- G_x;$ 
\end{enumerate}
If $G_x=e$ for every $x \in G$ we say that $G$ acts   {\em freely}. 
Note that  condition (1) implies that the quotient space $X/G$ is Hausdorff.  
It is well known that if $G$ acts properly discontinuously on $X$, 
then $X/G$ is a manifold and the projection $\pi:X \lra X/G$ 
is a covering map. 

In this article we study which groups act properly discontinuously on the
unit sphere of infinite dimensional Hilbert space and we 
shall prove that the situation is quite different from finite 
dimensional case. The main result (Theorem \ref{main}) 
is that any group of the form
$$H \cong G \oplus \Z_{p_1^{\alpha_1}} \oplus \cdots \oplus
\Z_{p_m^{\alpha_m}} ,$$
where $G$ is  a group without elements of finite order, 
acts isometrically and 
properly discontinuously on the unit sphere ${\bf S} (l_2(G))$ of $l_2(G)$
if and only if $p_i$ are different primes. Moreover, we shall prove    
that a group $H$ acts isometrically and properly discontinuously 
on the unit sphere of some separable Hilbert space 
then $H$  acts, with the same properties, on the unit sphere of any Hilbert 
space.
(Proposition \ref{aumenta}).
Some basic  references for
infinite dimensional geometry are \cite{La} and \cite{Kl}. 
\section{Properly discontinuous isometric actions on the infinite dimensional
unit sphere of a Hilbert space}   
Let $f$ be a mapping from a set $G$ into a real (complex) Banach space $\E$.
We recall, briefly, the idea of unordered summation in terms of 
convergence of nets. The set $\Lambda$ of finite subsets of $G$ is 
directed by the inclusion relation $\subseteq$, and we can define a net 
by the equation 
\[
\eta : \Lambda \lra \E, \ \ \eta(F)= \sum_{g \in F} f(g) . 
\]
We say that $f$ is summable if the net $\eta$  converges to some element, 
and this is possible if and only if there exists an element 
$L_o \in \E$ with the following property: for every 
$\epsilon > 0$ there exists 
$J(\epsilon) \in \Lambda$ such that
\[
\pr L_o \ - \ \sgj f(g) \pr < \epsilon, \ {\rm whenever} \ 
J(\epsilon) \subseteq J \in \Lambda.
\] 
It is well known, see \cite{Ru}, that 
$$l_2(G) = \{x: G \vai \R : \sum_{g\in G}[x(g)]^2 < \infty \},$$
is a Hilbert space with inner product 
$\langle x,y \rangle = \sg x(g)y(g)$ and a Hilbert 
basis is given  by the functions $e_h(g) = \delta_{hg}, h, g \in G$. 
When $G$ is countable then $l_2 (G)$ is isometrically isomorphic to $l_2$ and 
any Hilbert space is isometrically isomorphic to $l_2(D)$, 
where $D$ is a set with the same cardinality of any  Hilbert basis. 
Generally, we have the following properties:  
\begin{enumerate}
\item $N = \{g\in G : x(g) \neq 0 \}$ is countable;
\item $\sum_{g\in G}[x(g)]^2 = sup\{ \sum_{g\in F}[x(g)]^2:
F \in \Lambda \}$. 
\end{enumerate}
The second  property proves that for every bijective map  $\phi:G \lra G$, if 
$x \in l_2(G)$ then $x\ \circ \phi \in l_2 (G)$ and 
$\pr x \pr = \pr x \circ  \phi \pr$. 

Now, we assume that $G$ is a group and we will denote by $R_g$ the right
translation, i.e. $h \lra hg$, that is a bijective map.  
It is easy to prove that the following application 
$$\mu : G \times l_2(G) \vai l_2(G)\ \ \mu(g,x) =gx= x\circ R_g$$
is an isometric effective action of $G$ into $l_2 (G),$ 
i.e. if $gx=x$ for every $x$ then $g=e,$  and, in particularly, 
we have obtained an isometric action on the unit sphere 
${\bf S} (l_2 (G))$ of $l_2 (G)$. 
Summing up we have the following result.
\begin{prop} \label{iso}
Any group $G$ is isomorphic to a subgroup of the unitary group of the Hilbert
space $l_2(G)$.
\end{prop}
Now we shall prove the first step of our  main result.
\begin{thm} \label{primo}
If  $G$ has no elements of finite order, then the action $\mu$ 
on the unit sphere  ${\bf S} (l_2(G))$ is  properly  discontinuous. Moreover 
${\bf S} (l_2(G))/G$ has a complete Riemannian metric with constant sectional 
curvature $1$.    
\end{thm}
\n
\dem $\ $ we will prove the above Theorem in three simples steps. \\
{\em Step 1: the action is free.} 

If $gx = x$ then $g^nx = x$  and we have, for every 
$h \in G$  and $n \in \Z$, $x(h) = x(hg^n)$. 
The elements $\{ hg^n , n \in \Z \}$\ are different because
  $g$\ has not finite order. Then there exists an element   $h \in G$ such that
$x(h) \neq 0$ and we have 
$\sum_{n=1}^{\infty}[x(hg^n)]^2 = \sum_{n=1}^{\infty}[x(h)]^2 = \infty$,
which is absurd by property 2.\\
\n
{\em Step 2: if the action is not properly
discontinuous, there exists a  sequence  of different  
elements $(g_n) \subseteq G$ and 
$x \in {\bf S} (l_2(G))$ such that  $g_nx \vai y$, for some  $y$}.

Since $G$ acts isometrically, this is an easy consequence 
of the  definition of properly  discontinuous action. \\
\n
{\em Step 3: there are no  sequences
$(g_n) \subseteq G$ of different 
elements and $x \in {\bf S} (l_2(G))$ such that $g_nx$  converges to 
some  $y \in {\bf S} (l_2(G))$}.

Assume that there exists  a sequence $g_n \subseteq G$  and
$x \in {\bf S} (l_2(G))$ such that $g_n x \vai y$, for some $y \in S(L_2(G))$. 
Hence for every $\alpha \in G$ we have
$$y(\alpha) = \langle y,e_{\alpha} \rangle = lim_{n \vai \infty}
\langle g_nx,e_{\alpha} \rangle
= lim_{n \vai \infty} x(\alpha g_n)$$
and by property 2 we conclude that $y(\alpha)=0$, for every $\alpha \in G$, 
which is a contradiction.

Now, it is well known that $G$ induces a metric on
the manifold ${\bf S} (l_2(G))/G$ such that  it is an  Hilbert manifold 
and the projection 
$ \pi: {\bf S} (l_2(G)) \lra {\bf S} (l_2 (G)) /G$ is a local isometry and a  
covering map. Obviously the sectional curvature of 
${\bf S} ( l_2 (G))/G$ is constant
and equal to $1,$ so to conclude the proof we shall prove that 
${\bf S} (l_2(G))/G$ is a complete Hilbert  manifold. 
In the finite dimensional 
geometry we can prove it easily, 
because being geodesically complete at some point $p,$ i.e. the
$\exp_p$ being defined on $T_p M,$  
implies completeness as metric space.  Unfortunately, this fact doesn't 
hold in infinite dimensional
geometry: Atkin, see \cite{At}, gave an example of a complete Hilbert 
manifold $M$, such that the exponential map is not 
surjective in some point $p \in M$. 
Take $ q \not\in \exp_p (T_p M);$ clearly $M- \{ q \}$ is already geodesically 
complete in  $p$ but it is not complete as metric space.
In this case we will use a simple criterion to resolve our problem. 

We say  that a continuous curve $c:[a,b) \lra M$ is {\em convergent} 
if there exists a discrete subset $D \subset [a,b)$ such that
\begin{itemize}
\item   $c:[a,b)- D \lra M$ is differentiable;
\item   $\lim_{t \rightarrow b} \int_{a}^{t} \pr \fis c(t) \pr dt$ 
is finite.
\end{itemize}
\n
{\em A Hilbert manifold  is complete if and only if the trace of 
any convergent curve is relatively compact.} 

If $(M,g)$  is complete then it is easy to prove that any trace
of convergent curve is relatively compact. Vice-versa, we assume that any 
trace of convergent curve is relatively compact. 
We shall prove that any Cauchy 
sequence has a convergent subsequence, which proves that
$M$ is complete as metric 
space. Let $( x_n )$ be any Cauchy sequence. For any 
$\epsilon=\frac{1}{2^k}$ there exists
$n(k)$ such that $\forall n,m \geq 
n(k),\ \Rightarrow d(x_n,x_m) \leq \frac{1}{2^{k+1}}.$
Note that we shall assume $n(k+1) \geq n(k)\ +\ 2$. Let  $( a_n )$ be a 
sequence such that
\[
\left \{ \begin{array}{l}
a_1=0 ; \\
a_n > 0,\ n \geq 2 ; \\
\sum_{1}^{\infty} a_n =1.
\end{array}
\right.
\]
By definition of distance in a Hilbert manifold, 
for any $k \in \N$ there exists a differential curve,
$\gamma_k:[\sum_{1}^{k} a_n ,\sum_{1}^{k+1} a_n ] \lra M$, between
$x_{n(k)}$ and  $x_{n(k+1)}$ such that  
\[
L[\gamma_k] \leq d(x_{n(k)},x_{n(k+1)})\ + \ \frac{1}{2^{k+1}}
\leq \frac{1}{2^k}.
\]
We define $\gamma:[0,1) \lra M,$ 
$\gamma(t)=\gamma_k (t)$  if  $t\in [\sum_{1}^{k}
a_n,\sum_{1}^{k+1} a_n]$. It is easy to see that 
$\gamma$ is a convergent curve so $( x_n )$ has a convergent
subsequence.  

Now take a convergent curve $ \gamma;[a,b] \lra M$ in ${\bf S} (l_2 (G))/G.$ 
We can lift $\gamma$ to a curve $\overline{\gamma}$ in 
$ {\bf S} (l_2(G))$ that is
a convergent curve because $\pi$ is a local isometry. On the other hand 
$ \gamma ([a,b]) \subseteq  \pi (\overline{ \overline{\gamma} ([a,b]) } ),$ 
then the trace of $\gamma$ is relative compact. Q.E.D.

Now,  the  infinite dimensional sphere is contractible, 
by Bes\-sega Theorem \cite{Be}, 
so we have proved the following result. 
\begin{cor}
Let $G$ be a group without elements of finite order. 
Then there exists a complete Hilbert
manifold with positive constant sectional curvature with fundamental group $G$.
\end{cor}
The above  result is different from the corresponding one in finite dimensions:
by Bonnet Theorem the fundamental group of a complete finite dimensional 
manifold with constant positive curvature is finite. 
%
%
On the contrary,  groups as $\Q$ and $\Z$ act  
properly discontinuously and isometrically on the unit sphere of $l_2$.

An interesting fact is that, when the model space is separable, then the 
cardinality  of the fundamental group is at most countable. 
This is a trivial consequence  of the lemma of Sierpinski \cite{Si}, 
which says that
every connected locally separable metric space  is separable, and of 
the fact that the cardinality of the fiber of the universal covering map, is 
the cardinality of the fundamental group. In particular our Theorem is
sharp, relative to the cardinality of the fundamental group, in a case of
Hilbert manifold modeled on $l_2.$

Now, we shall study the case when the group $G$ has  elements of finite order 
and  we shall investigate what are the necessary conditions for a 
linear map,  with finite order, to have no fixed points except for the origin. 
In the finite dimensional case, the 
groups  acting  properly discontinuously, as isometric groups, on the unit 
sphere are classified \cite{Wo}  and the first question is: if a  group 
acts on the finite dimensional sphere, does it act on any  
infinite dimensional Hilbert space?
The following proposition gives a positive answer to this question.
\begin{prop} \label{aumenta}
If $G$ acts isometrically and properly discontinuously on the unit sphere of 
$\R^{n}$ or $l_2$  then $G$ 
acts  with the same properties on the unit sphere of any infinite dimensional 
Hilbert space.
\end{prop}
\n
\dem $\ $ first of all we shall prove that if  
$G \subset \bf{O} (n)$,  then it acts on 
the unit sphere of $l_2$ with the same properties. We define
\[
T_g : l_2 \lra l_2,\ \ \ \ T_g
(\sum_{i=1}^{\infty} x_i e_i) =
\sum_{i=1}^{\infty} g (\sum_{j=n(i-1)}^{ni-1} x_{j+1} e_{j+1}),
\]
with the following identification:
\[
R^n \stackrel{T_i}{\lra} l_2, \ \ \ 
(x_1,\ldots,x_n) \lra \sum_{j=1}^{n} x_{j} e_{n(i-1)+k-1}.
\]
Clearly $T_g$ is an isometry and it is easy to check that this 
defines an isometric action on ${\bf S} (l_2)$, 
that is properly discontinuous since $G$ is finite. 
The general case is analogous. We know that any
Hilbert space $\H$ is isometrically isomorphic to
 $l_2(D)$, where $D$ is a set that  has the same cardinality of any
Hilbert basis of $\H$, and 
\[
D=\bigcup_{i \in I} D_i,\ D_i \bigcap D_j = \emptyset,\ {\rm if}
\ i \neq j, \ D_i \cong  \N
\]
(see \cite{La2} page 678). Let $f \in l_2 (D)$. We define
\[
f_i (x) = \left \{
\begin{array}{ll}
f(x) & {\rm if}\  x \in D_i; \\
0 & \textrm{otherwise.} \\
\end{array}
\right.
\]
Now, $f \in l_2(D)$ and by property $2$ we have that $f_i \in l_2 (D)$.
Another consequence of property $2$, that is easy to check, is that 
for every  $\epsilon>0$
there exist a $J(\epsilon) \in \Lambda$ such that  for every finite set
$J \subseteq D- J(\epsilon)$ we have  
$\sum_{d \in J} [f(d)]^2 \leq \epsilon$. Let $\epsilon>0$ and let 
$J(\epsilon)$ be above. Let 
$I_o:= \{ i \in I:\ D_i \bigcap J(\epsilon) \neq \emptyset \}$ and we compute,
for every finite set $I \supseteq I_o$,
\begin{small}
$$
\begin{array}{lcl}
\pr \sum_{i \in I} f_i \ - \ f \pr^2 & = & sup \{ \sum_{g \in J} [f_i (g) \ -
f(g)]^2, \ {\rm J \ finite} \} \\
                    &=& sup \{ \sum_{g \in J- J(\epsilon)} [f(g)]^2, \ {\rm J
\ finite} \} \\
                    &\leq & \epsilon.
\end{array}$$
\end{small}
Hence $\sum_{i \in I} f_i$ converges to $f$ and it is easy to check that
$\alpha f \ + \beta g =\sum_{i \in I} \alpha f_i \ + \ \beta g_i$ and
$f=\sum_{i \in I} f'_i$, where $f'_i=0$ on the complement of 
$D_i$, so  $f_i=f'_i$. Moreover, for every finite set $J \subseteq I$ 
we have 
\begin{small}
$$\begin{array}{lcl}
\sum_{i \in J} \pr f_i \pr^2 & = & \sum_{i \in J} sup (
\sum_{ g_i \in  J_i} [f_i (g_i)]^2) \ \ \ J_i \subseteq D_i \ {\rm finite} \\
                    &=& sup (\sum_{g \in \cup J_i} [f(g)]^2)
\end{array}$$
\end{small}
so $\sum_{i \in I} \pr f_i \pr^2 = \pr f \pr^2$. Now, if we define
\[
T_g : l_2 (D) \lra l_2(D), \ \
T_g (\sum_{i \in I} f_i)=\sum_{i \in I} g f_i,
\]
by the identifications  $l_2 \stackrel{T_i}{\lra} l_2 (D)$ as
\[
f \lra \left \{ 
\begin{array}{ll}
f \circ s_i (x) &  \textrm{if} \ x \in D_i; \\
0               &   \textrm{otherwise,} \\
\end{array}
\right.
\]
where $s_i$ is a bijective map between $D_i$ and $\N$ , we have an isometric 
action that is free.
If $G$ does not act 
properly discontinuously, there will be a sequence  $g_n$ of distinct 
elements and an  element
$x=\sum_{i \in I} x_i \neq 0$ such that
\[
g_n x=\sum_{i \in I} g_n x_i \vai y=\sum_{i \in I} y_i.
\]
Hence $g_n x_i \vai y_i$, that is a contradiction.  Q.E.D. 

Let $\H$ be a Hilbert space and let $T: \H \lra \H$ be a linear continuous
map with finite order $m,$ i. e. $T^m=Id$ and $m$ is the smallest integer with 
this property. Henceforth, when we say that a liner map has no fixed 
points, we will mean that the unique fixed point is the origin. If $T$ has no 
fixed point, then it must  satisfy the following equation
\[
T^{m-1} + T^{m-2} + \cdots +T + Id = 0.
\] 
If  $m = 2$ then $T=-Id.$ Hence, if a group $H$ acts without fixed point then
it has at most one  element $g$ of order two and $g$ has to  belong 
to the  center of $H$. 
In particular, groups as   $\Z_2 \oplus \Z_2$,
the symmetric group 
$ S_3 = <\phi , \psi ; \phi ^2 = \psi^3 = 1, \phi \psi = \psi ^2 \phi >$ 
and the diedral group
$ D_8 = <\sigma , \tau ; \sigma ^4 = \tau^2, \tau\sigma\tau = \sigma ^3>,$ 
cannot act  without fixed point.

We shall assume that $\H$  is a complex Hilbert space. Otherwise we shall 
consider the complexification $\H\otimes \C$ and $T^{\C}$  and, clearly, 
$T$ has fixed points if and only if  $T^{\C}$ does. 
Since $T^m =Id,$ the minimal polynomial is well-defined and is given by
$p_T (z)= (z - \xi_1)\cdots(z - \xi_n),$ where
 $\xi_i$ are roots of unity. As in the finite dimensional case,
we have the following decomposition 
$$\H = \H_1 \oplus \cdots \oplus \H_n$$
where 
\begin{itemize}
\item $\H_i = Ker(T - \xi_iId)$;
\item $p_{T|\H_i}(z) = (z - \xi_i)$.
\end{itemize}
Moreover, 
 $x =x_1\ + \cdots + x_n$ is a fixed point of  $T^k$, $1 \leq k \leq m$ if
and only if $T^k (x_i)=x_i$. This simple remark proves the following result.
\begin{prop}
Let $T:\H \vai \H$ be a linear continuous map with finite order $m$. 
Then $T, \cdots, T^{m-1}$ have no fixed points if and only if 
\[
\Phi_m(T)=0
\]
where $\Phi_m(t)$ is the $m$-th cyclotomic polynomial.
\end{prop}

Now, let $T:\H \lra \H$ be a linear continuous map of order $n$ and let  
$L:\H \vai \H$ be a linear continuous map of order $m$. 
We suppose that $L \circ T= T \circ L$ and  
there exists a prime  $p$ such that  $p|m$ and  $p|n$. 
Let $T_1=T^{\frac{n}{p}}$, 
$L_1=L^{\frac{m}{p}}$. Clearly the order of the last two linear operators is 
$p$, and $L_1 \circ T_1=T_1 \circ L_1$. Relative to $T_1$ we have 
\[
\H=\H_1 \oplus \cdots \oplus \H_l,
\]
and $H_i$ are $L_1$-invariant. In particular, relative to $L_1$ we have
\[
\H_1=V_1 \oplus \cdots \oplus V_r.
\]
Now it is easy to prove the following simple, but important, result.
\begin{prop} \label{prop}
Let $G$ be a group that acts on $\H$ without fixed points, except the origin, 
as linear group. Let $g \in G$ be an element of finite order and 
let $m \in C(g)$ (centralizer of $g$) have  finite order.
Then, either $g^s = m^k$, for some $k,s \in \N$ or $MDC(o(m),o(g))=1$.
\end{prop}
The above result proves that any group that  has a subgroup of the form 
\[
\Z_n \oplus \Z_m
\]
with  $MDC(n,m) \neq 1$,  does not act, as linear group, 
without fixed points.

Now, we shall construct explicitly an action of the group 
\[
 \Z_{p_1^{\alpha_1}} \oplus \cdots \oplus
\Z_{p_m^{\alpha_m}}
\]
with  $p_1< \ldots <p_m$ primes. Let $\xi_1 \ldots \xi_m$ be primitive
$p_i^{\alpha_i}$-th roots of unity and $\theta_i=arg(\xi_i)$. 
Generally for any $\theta$ we define  
\begin{small}
\[
T_{\theta}:\H \oplus \H \lra \H \oplus \H,\ \ 
T_{\theta}(x,y)=(\cos (\theta) x \ - \ \sin (\theta)y,
         \sin (\theta) x \ + \ \cos (\theta)y).
\]
\end{small}
It is easy to check that
\begin{enumerate}
\item $T_{\theta}$ is a surjective isometry;
\item $\forall \theta, \alpha \in \R$,
$T_{\theta} \circ T_{\alpha}= T_{\alpha} \circ T_{\theta}=
T_{\theta + \alpha}$;
\item $T_{\theta}$ has fixed points if and only if $\theta= 2k \pi$, with
$k \in \Z$.
\item $\prod_{i=1}^{m} \xi_i^{s_i} (x\ + \ iy)=
\pi_1(\prod_{i=1}^{m} T_{\theta_i}^{s_i} (x,y)) \ + \
i(\pi_2(\prod_{i=1}^{m} T_{\theta_i}^{s_i} (x,y)))$,
 where $\pi_1$ and $\pi_2$  are the natural projections.
\end{enumerate} 
We shall prove that the operators $T_{\theta_i}$ represent     
$\Z_{p_1^{\alpha_1}} \oplus \cdots \oplus \Z_{p_m^{\alpha_m}}$ 
as isometry subgroup of $\bf{O} ( \H \oplus \H)$.
Suppose there exist  $s_1, \ldots, s_m \in \N$ 
such that
\[
T_{\theta_i}^{s_i}= \prod_{j \neq i} T_{\theta_j}^{s_j}.
\]
Then by the  property 4 we have
\[
\xi_i^{s_i} = \prod_{j \neq i} \xi_j^{s_j}
\]
which implies that $s_1= \ldots =s_m=0$, because $p_i$ are different primes. 
If
\[
\prod^{m}_{i=1}T_{\theta_i}^{s_i}
\]
has fixed points then we have, by property 4, 
\[
\prod^{m}_{i=1} \xi_i^{s_i}=1
\]
that implies, by the same  argument as above, $s_1= \ldots =s_m=0$.
Now we shall prove the principal result.
\begin{thm} \label{main}
Let $H \cong G \oplus \Z_{p_1^{\alpha_1}} \oplus \cdots \oplus
\Z_{p_m^{\alpha_m}}$ be a group where $G$ is  a group without elements of 
finite order. Then  $H$ acts isometrically and 
properly discontinuously  on the unit sphere ${\bf S} (l_2(G))$ of $l_2(G)$
if and only if $p_i$ are different primes. Furthermore, 
${\bf S} (l_2 (G))/H$ has a complete Riemannian metric with constant sectional
curvature 1.  
\end{thm} 
\n
\dem $\ $ the forward implication  has already been verified 
in Proposition \ref{prop}. For the reverse, we shall prove that  
$H$ acts isometrically on the Hilbert space
$l_2(G) \oplus l_2(G)$ that is isometric to  
 $l_2(G)$. Any element  $g \in G$ acts on  $l_2(G) \oplus l_2(G)$ as 
$g \oplus g$, and any $g$ acts as in Theorem \ref{primo}. 
It is easy to check 
that this action has the same properties of the action 
of $G$ into  $l_2(G)$. The elements of finite order, 
that we will indicated as $A$,
acts as above. The proof is divided in three steps.
\begin{enumerate}
\item $G \bigcap A=e$. \\
If $\prod^{n}_{i=1} T_{\theta_i}^{s_i}=T_{\theta}=g_o$, where
$\theta= \theta_1 s_1 \ + \cdots \ + \theta_m s_m$, then
\[
\begin{array}{lcl}
g_o(e_e,e_e) & = & (e_{{g_o}^{-1}},e_{{g_o}^{-1}}) \\
             & = & (e_e (\cos \theta \ - \sin \theta),
                   e_e (\cos \theta \ + \sin \theta)). \\
\end{array}
\]
Then $g_o=e$ and $\theta= 2 \pi k$;
\item $\forall g \in G$ e $\forall a \in A$, $ag=ga$.\\
We shall prove that  
$T_{\theta}$ commutes  with  $g$, for every $g \in G$ and every 
 $\theta \in \R$.\\
\begin{small}
$g(T_{\theta} ((\sum_{h \in G} x_h e_h, \sum_{h \in G} y_h e_h)) =$\\
$g(\sum_{h \in G} (x_h \cos \theta \ - \ y_h \sin \theta) e_{h},
\sum_{h \in G} x_h \sin \theta \ + \ \cos \theta) e_{h}))
=$ \\
$
(\sum_{h \in G} x_h\cos \theta \ - \ y_h \sin \theta) e_{hg^{-1}},
\sum_{h \in G} (x_h \sin \theta \ + \ \cos \theta) e_{hg^{-1}}))= $ \\
$
 T_{\theta}((\sum_{h \in G} x_h e_{hg^{-1}}, \sum_{h \in G}
y_h e_{hg^{-1}}))= $ \\
$T_{\theta}(g(\sum_{h \in G} x_h e_h, \sum_{h \in G} y_h e_h))).
$
\end{small}
\item $ag$  has no fixed points and the action is properly discontinuous. \\
If
\[
ag(\sum_{h \in G} x_h e_h, \sum_{h \in G} y_h e_h)=
(\sum_{h \in G} x_h e_h, \sum_{h \in G} y_h e_h)
\]
then
\[
\sum_{h \in G} (x_h\cos \theta \ - \ y_h \sin \theta)e_{hg^{-1}}=
\sum_{h \in G} x_h e_h
\]
and
\[
\sum_{h \in G} (x_h\cos \theta \ + \ y_h \sin \theta)e_{hg^{-1}}=
\sum_{h \in G}  y_h e_h.
\]
Hence, for every $\alpha \in G$ we have 
\[
x_{\alpha} \cos \theta \ + y_{\alpha} \cos \theta =y_{\alpha g^{-1}}; \ \ (1)
\]
\[
x_{\alpha} \cos \theta \ - y_{\alpha} \cos \theta =x_{\alpha g^{-1}}. \ \ (2)
\]
Now, there exists an  $\alpha \in G$
such that  $x_{\alpha} \neq 0$ or $y_{\alpha} \neq 0$. Then, it easy to check 
by (1)  and (2),  that for every $n \in \N$ we have
\[
| x_{\alpha g^{-n}} |^2 \ + \
| y_{\alpha g^{-n}} |^2 \ =| x_{\alpha}|^2 \ + \ | y_{\alpha} |^2.
\]
Hence 
\[
\sum^{\infty}_{n=1}  | x_{\alpha g^{-n}} |^2
\ + \
\sum^{\infty}_{n=1}| y_{\alpha g^{-n}} |^2 = \infty          
\]
which is a contradiction, because the above series have to converge.
If the action is not properly discontinuous there exists a sequence 
of distinct elements $g_n a_n$  and an element $x \neq 0$,
such that $g_n a_n x$  converge to some element. The group $A$ is finite, so
we shall consider a subsequence such that 
$a_{n(k)}=a_o$. Then we have obtained a  contradiction because $G$ acts 
properly discontinuously on the unit sphere. Now, we conclude our proof as
in Theorem 1.  Q.E.D.
\end{enumerate}

A trivial application of main Theorem  and Proposition \ref{aumenta}
is the following
result that gives a necessary and sufficient conditions when the group $G$
is a finitely generated Abelian group.
\begin{cor} \label{cormain}
Every finitely generated Abelian group $G$ acts isome\-tri\-cal\-ly 
and properly 
discontinuously on the unit sphere of any infinite dimensional Hilbert space 
if and only if $G$ has a torsion of the form
\[
\Z_{p_1^{\alpha_1}} \oplus \cdots \oplus \Z_{p_m^{\alpha_m}},
\]
where $p_i$ are distinct primes. In particular,  there exists a 
complete Hilbert manifold, possibly modeled on $l_2$,
with  positive constant sectional curvature, whose fundamental group is $G$.
\end{cor}

$\ $ \\
$\ $ \\
Leonardo Biliotti \\
Dipartimento di Matematica e delle Applicazioni all'Architettura, \\
Universit\`a degli studi di Firenze, \\
Piazza Ghiberti 27 - Via dell'Agnolo 2r - 50132 Firenze (Italy) \\
e-mail:{\tt biliotti@ime.unicamp.br} \\
\end{document}